\documentclass[12pt,verbatim]{amsart}
\usepackage[pdftex]{graphicx}
\usepackage{amsfonts}
\usepackage{amsthm}
\usepackage{amsmath}
\usepackage{amssymb}
\usepackage{amscd}

\newcounter{minutes}
\divide\time by 60
\newcounter{hours}
\multiply\time by 60 \addtocounter{minutes}{-\time}
\begin{document}


\title{}{}
\author{}{}

\centerline{\large \bf Numerical methods with Sage}
\bigskip
\centerline{ \bf Lauri Ruotsalainen and
Matti Vuorinen}

\email{lauri.ruotsalainen@gmail.com, vuorinen@utu.fi}

\address{Department of Mathematics and Statistics, University of Turku}


\maketitle

\def\thefootnote{}
\footnotetext{ \texttt{\tiny File:~\jobname .tex,
          printed: \number\year-\number\month-\number\day,
          \thehours.\ifnum\theminutes<10{0}\fi\theminutes }
} \makeatletter\def\thefootnote{\@arabic\c@footnote}\makeatother

\begin{abstract}
Numpy and SciPy are program libraries for the Python scripting language, which apply to a large spectrum
of numerical and scientific computing tasks. The Sage project provides a multiplatform software environment
which enables one to use, in a unified way, a large number of software components, including Numpy and Scipy, and which has Python as its command language. We review several examples, typical for scientific computation courses, and their solution using these tools in the Sage environment.
\end{abstract}


\section{Introduction}

Python is a popular multipurpose programming language. When augmented with special libraries/modules, it is suitable also for scientific computing and visualization. These modules make scientific computing with Python similar to what commercial packages such as MATLAB, Mathematica and Maple offer. Special Python libraries like NumPy, SciPy and CVXOPT allow fast machine precision floating point computations and provide
subroutines to all basic numerical analysis functions. As such the Python language can be used to combine several software components. For data visualization there are numerous Python based libraries such as Matplotlib. There is a very wide collection of external libraries. These packages/libraries are available either for free or under permissive open source licenses which makes these very attractive for university level teaching purposes.

The mission of the Sage project is to bring a large number of program libraries under the principle of the GNU General Public License under the same umbrella. Sage offers a common interface for all its components and its command language is Python. Originally developed for the purposes of number theory and algebraic geometry, it presently supports more than 100 special libraries. Augmented with these, Sage is able to carry out both symbolic and numeric computing. \cite{Sage}

In this paper, we present examples or case studies of the usage of Sage to solve some problems common for typical numerical analysis courses. The examples are drawn from the courses of the second author at the University of Turku, covering the standard topics of numerical computing and based, to a large extent, on the standard textbooks \cite{BF,CdB,H,MF,Mol,K,L,TLN}.

\section{General observations}

\subsection{History of Sage}

Sage was initially created by William Stein in 2004-2005, using open source programs released under the GPL or a GPL-compatible license. The main goal of the project was to create a viable open source alternative to proprietary mathematical software to be used for research and teaching. The first release of the program was in February 2005. By the end of the year, Sage included Pari, GAP and Singular libraries as standard and worked as a common interface to other mathematical programs like Mathematica and Magma. \cite{His, PARI, GAP, DGPS}

Since the beginning Sage has expanded rapidly. As of July 2012, at least 246 people have actively contributed code for Sage. The range of functionality of Sage is vast, covering mathematical topics of all kinds ranging from number theory and algebra to geometry and numerical computation. Sage is most commonly used in university research and teaching. On Sage's homepage there are listed more than one hundred academic articles, books and theses in which the program has been involved.

\subsection{Access to Sage}

Sage can be utilized in multiple ways. The main ways to use Sage are:

\begin{enumerate}

\item Over the network. In this case no installation is required.

\item The software is downloaded and installed on the personal work station.

\end{enumerate}

\noindent These ways to use Sage will now be described more closely.

\begin{enumerate}

\item One of the strengths of Sage is that it can be used over the network without requiring any installation of the application. The Sage Notebook is a web browser-based graphical user interface for Sage. It allows writing and running code, displaying embedded two and three dimensional plots, and organizing and sharing data with other users. The Sage Notebook works with most web browsers without the need for additional add-ons or extensions. However, Java is required to run Jmol, the Java applet used in Sage to plot 3D objects.

\item Sage can be also installed natively for Linux, OS X and Solaris. Both binaries and source code are available for download on Sage's homepage. In order to run Sage on the Windows operating system the use of virtualization technology (e.g. VirtualBox or VMware) is required. There are three basic interfaces to access Sage locally: the Sage Notebook on a web browser, a text-based command-line interface using IPython, or as a library in a Python program. \cite{PG}

\end{enumerate}

\noindent In addition to these main options, Sage can be applied in alternative ways. For instance, a single cell containing Sage code can be embedded in any webpage and evaluated using a public single cell server. There is also support to embed Sage code and graphics in LaTeX documents using Sagetex. \cite{SageTeX}

\subsection{Key issues}

Numerical computation has been one of the most central applications of the computer since its invention. In modern society, the significance of the speed and effectiveness of the computational algorithms has only increased with applications in data analysis, information retrieval, optimization and simulation.

Most numerical algorithms are implemented in a variety of programming languages. There are various commercial collections of software libraries with numerical analysis functionality, mostly written in Fortran, C and C++. These include, among others, the IMSL and NAG libraries. Computer algebra systems and applications developed specially for numerical computing include MATLAB, S-PLUS, Mathematica, IDL and LabVIEW. There are also free and open source alternatives such as the GNU Scientific Library (GSL), R, Scilab, Freemat and Sage.

Many computer algebra systems (e.g. Maple, Mathematica and MATLAB) contain an implementation of a new programming language specific to that system. In contrast, Sage uses Python, which is a popular and widespread high-level programming language. The Python language is considered simple and easy to learn while making the use of more advanced programming techniques, such as object-oriented programming and defining new methods and data types, possible in the study of the mathematics. Python functions as a common user interface to Sage's nearly one hundred software packages.

Some of the advantages of Sage in scientific programming are the free availability of the source code and openness of development. Most commercial software do not have their algorithms publicly available, which makes it impossible to revise and audit the functionality of the code. Therefore the use of the built-in functions of these programs may be inadequate in some mathematical studies that are based on the results given by these algorithms. This is not the case in open source software, where the user can verify the individual implementations of the algorithms.

In many situations, the Python interpreter is fast enough for common calculations. However, sometimes considerable speed is needed in numerical computations. Sage supports Cython, which is a compiled language based on Python that supports calling C functions and declaring C types on variables and class attributes. It is used for writing fast Python extension modules and interfacing Python with C libraries. Significant parts of the NumPy and SciPy libraries included with Sage are written in Cython. In addition, techniques like distributed computing and parallel processing using multi-core processors or multiple processors are supported in Sage. \cite{BBSE, Num}

Sage includes the Matplotlib library for plotting two-dimensional objects. For three-dimensional plotting, Sage has Jmol, an open-source Java viewer. Additionally, the Tachyon 3D raytracer may be used for more detailed plotting. Sage's \textsl{interact} function can also bring added value to the study of numerical methods by introducing controllers that allow dynamical adjusting of the parameters in a Sage program. \cite{Hun, Jmol, T, Tac}

\subsection{Numerical tools in Sage}
\label{sec:numerical tools}

Sage contains several program libraries suitable for numerical computing. The most substantial of these are NumPy, SciPy and CVXOPT, all of which are extension modules to the Python programming language, designed for specific mathematical operations. In order to use these packages in Sage, they must be first imported to the Sage session using the \textsl{import} statement. \cite{ADV, O, SciPy}

NumPy provides support for fast multi-dimensional arrays and numerous matrix operations. The syntax of NumPy resembles the syntax of MATLAB in many ways. NumPy includes subpackages for linear algebra, discrete Fourier transforms and random number generators, among others. The SciPy library is a library of scientific tools for Python. It uses the array object of the NumPy library as its basic data structure. SciPy contains various high level scientific modules for linear algebra, numerical integration, optimizing, image processing, ODE solvers and signal processing, et cetera.

CVXOPT is a library specialized in optimization. It extends the built-in Python objects with dense and sparse matrix object types. CVXOPT contains methods for both linear and nonlinear convex optimization. For statistical computing and graphics, Sage supports the R environment, which can be used via the Sage Notebook. Some statistical features are also provided in SciPy, but not as comprehensively as in R. \cite{R}

For more information on numerical computing with Sage, see \cite{Num}.


\section{Case studies}

In this section our goal is to present examples of the use of
Sage for numerical computing. This goal will be best achieved
by giving code snippets or programs that carry out some typical
numerical computation task. We cover some of the main
aspects of a standard first course in numerical computing, such
as the books of \cite{BF,CdB,H, MF,Mol,K,L,TLN}.

\subsection{Newton's method}

Computing the solution of some given equation is one of the fundamental problems of numerical analysis. If a real-valued function is differentiable and the derivative is known, then Newton's method may be used to find an approximation to its roots.

In Newton's method, an initial guess $x_0$ for a root of the differentiable function $f: \mathbb{R} \rightarrow \mathbb{R} $ is made. The consecutive iterative steps are defined by $$x_{k+1} = x_k-\frac{f(x_k)}{f'(x_k)}, k=0,1,2,\dots$$

An implementation of the Newton's method algorithm is presented in the next code. As an example, we use the algorithm to find the root of the equation $x^2-3=0$. The function \textsl{newton\_method} is used to generate a list of the iteration points. Sage contains a preparser that makes it possible to use certain mathematical constructs such as $f(x) = f$ that would be syntactically invalid in standard Python. In the program, the last iteration value is printed and the iteration steps are tabulated. The accuracy goal for the root $2h$ is reached
 when $f(x_n-h)f(x_n+h)<0\,.$ In order to avoid an infinite loop, the maximum number of iterations is limited by the parameter \textsl{maxn}.\\

\begin{figure}[h]
	\centering
		\includegraphics[width=0.70\textwidth]{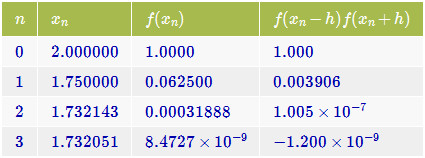}
	\caption{Output of the Newton iteration.}
	\label{fig:newton}
\end{figure}


\small
\begin{verbatim}
def newton_method(f, c, maxn, h):
    f(x) = f
    iterates = [c]
    j = 1
    while True:
        c = c - f(c)/derivative(f(x))(x=c)
        iterates.append(c)
        if f(c-h)*f(c+h) < 0 or j == maxn:
            break
        j += 1
    return c, iterates

f(x) = x^2-3
h = 10^-5
initial = 2.0
maxn = 10
z, iterates = newton_method(f, initial, maxn, h/2.0)
print "Root =", z
html.table(([i, c.n(digits=7), f(c).n(digits=5),
       (f(c-h)*f(c+h)).n(digits=4)] for i, c in enumerate(iterates)),
       header=["$n$", "$x_n$", "$f(x_n)$", "$f(x_n-h)f(x_n+h)$"])
\end{verbatim}
\normalsize

\subsection{Computational Methods of Linear Algebra}

Sage offers a good platform for practicing both symbolic and numerical linear algebra. The software packages specialized in computational linear algebra that are contained in Sage include LAPACK, LinBox, IML, ATLAS, BLAS and GSL. Both the NumPy and SciPy libraries have subpackages for algorithms that utilize NumPy's fast arrays for matrix operations. Also, the GSL library can be used via Cython. In most of the examples of this chapter, we use native Sage matrices that suit most needs.

Let $A$ be a matrix over the real double field (RDF) as defined in the next code. The \textsl{matrix} function accepts the base ring for the entries and the dimensions of the matrix as its parameters. We can compute various values associated with the matrix, such as the determinant, the rank and the Frobenius norm:

\begin{verbatim}
Sage: A = matrix(RDF, 3, [1,3,-3, -3,7,-3, -6,6,-2])
Sage: A.determinant()
  -32.0
Sage: A.rank()
  3
Sage: A.norm()
  12.7279220614
\end{verbatim}


The function A.inverse() returns the inverse of A, should it exist. Otherwise Sage informs that the matrix must be nonsingular in order to compute the inverse.

\begin{verbatim}
Sage: A.inverse()
\end{verbatim}

\small
$\left(\begin{array}{rrr}
-0.125 & 0.375 & -0.375 \\
-0.375 & 0.625 & -0.375 \\
-0.75 & 0.75 & -0.5 \\
\end{array}\right)\\$
\normalsize

Let $b=(1,3,6)$. We solve the matrix equation $Ax=b$ using the function \textsl{solve\_right}.  The notation $A\backslash b$, specific to Sage, may also be used.

\begin{verbatim}
Sage: b = vector([1,3,6])
Sage: A.solve_right(b)
  (-1.25, -0.75, -1.5)
\end{verbatim}

In numerical linear algebra, different decompositions are employed to implement efficient matrix algorithms. Sage provides several decomposition methods related to solving systems of linear equations (e.g. LU, QR, Cholesky and singular value decomposition) and decompositions based on eigenvalues and related concepts (e.g. Schur decomposition, Jordan form). The availability of these functions depends on the base ring of matrix; for numerical results the use of real double field (RDF) or complex double field (CDF) is required.

Let us determine the LU decomposition of the matrix $A$. The result of the function \textsl{A.LU()} is a triple of matrices $P$, $L$ and $U$, so that $PA=LU$, where $P$ is a permutation matrix, $L$ is a lower-triangular matrix and $U$ is an upper-triangular matrix.

\begin{verbatim}
Sage: A.LU()
\end{verbatim}
\small
$$\left(\left(\begin{array}{rrr}
0.0 & 0.0 & 1.0 \\
0.0 & 1.0 & 0.0 \\
1.0 & 0.0 & 0.0
\end{array}\right), \left(\begin{array}{rrr}
1.0 & 0.0 & 0.0 \\
0.5 & 1.0 & 0.0 \\
-0.1666 & 1.0 & 1.0
\end{array}\right), \left(\begin{array}{rrr}
-6.0 & 6.0 & -2.0 \\
0.0 & 4.0 & -2.0 \\
0.0 & 0.0 & -1.3333\\
\end{array}\right)\right)$$
\normalsize
\newline

According to linear algebra, the solution of the equation $Ax=b$ for a $n \times n$ matrix is unique if the determinant $det(A) \neq 0\,.$ However, the solution of the equation may be numerically unstable also if $det(A) \neq 0\,.$ The standard way to characterize the "numerical nature" of a square matrix $A$ is to use its condition number $cond(A)$ defined as $\sigma_M/\sigma_m$ where   $\sigma_M$   ($\sigma_m$) is the largest (least) singular value of  $A\,.$ The {\it singular value decomposition (SVD)} of $A$ yields the singular values as follows: $$   A= U S  V^T$$ where $U$ and $V$ are orthogonal $n \times n$ matrices and the $S$ is a diagonal $n \times n$ matrix with positive entries organized in decreasing order, the singular values of $A\,.$

In the next example we study the influence of the condition number on the accuracy of the numerical solution of a random matrix with a prescribed condition number. For this purpose we use a simple method to generate random matrices with a prescribed condition number $c\ge 1$:
take a random square matrix  $A$, form its SVD  $A= U S V^T $ and modify its singular values $S$ so that for the modified matrix $S_c$ the quotient of the largest and least singular value is $c$ and then $A_c= U S_c V^T$ is our desired random matrix with $cond(A_c) =c \,.$ For several values of $c$ we then observe the error
in the numerical solution of $A_c \, x =b$ and graph the error as a function of $cond(A_c)$ in the loglog scale. We
see that the condition number appears to depend on $cond(A_c)$ almost in a linear way.

\small
\begin{verbatim}
from numpy import *
from matplotlib.pyplot import *
data = []
n = 20
A = random.rand(n,n)
U, s, V = linalg.svd(A)
ss = zeros((n,n))

for p in arange(1,16,2.):
    c = 10.^p
    for j in range(n):
        ss[j, j] = s[0] - j * (s[0] - s[0] / c) / (n - 1)
    aa = dot(dot(U, ss), V.T)
    b = dot(aa, ones(n))
    numsol = linalg.solve(aa, b)
    d = linalg.norm(numsol-ones(n))
    data.append([c, d])

data = array(data)
x,y = data[:,0],data[:,1]

clf()
loglog(x, y, color='k', linewidth=2)
loglog(x, y, 'o', color='k', linewidth=10)
xlabel('Condition number of the matrix',\
fontweight='bold', fontsize=14)
ylabel('Error', fontweight='bold', fontsize=14)
title('Error as a function of the condition matrix',\
fontweight='bold', fontsize=14)
grid(True)
savefig("fig.png")
\end{verbatim}
\normalsize


\begin{figure}[h!t]
	\centering
		\includegraphics[width=1.00\textwidth]{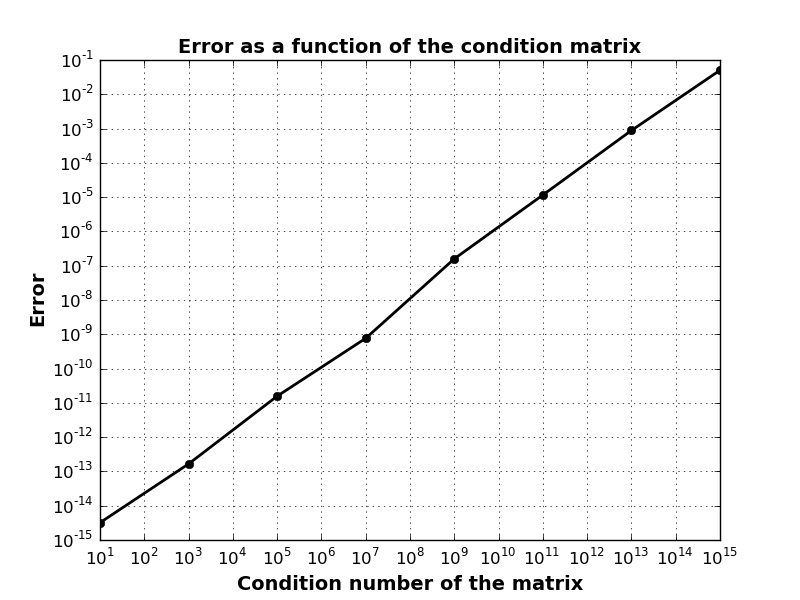}
	\caption{The output of the program used to study the influence of the condition number on the accuracy of the numerical solution of a random matrix with a prescribed condition number.}
	\label{fig:svd}
\end{figure}

\subsection{Numerical integration}

\begin{figure}[h!]
	\centering
		\includegraphics[width=1.00\textwidth]{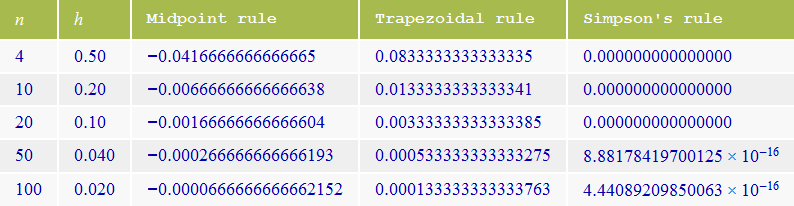}
	\caption{The table shows the difference between the exact value of the integral and the approximation using various rules.}
	\label{fig:integration}
\end{figure}

Numerical integration methods can prove useful if the integrand is known only at certain points or the antiderivate is very difficult or even impossible to find. In education, calculating numerical methods by hand may be useful in some cases, but computer programs are usually better suited in finding patterns and comparing different methods. In the next example, three numerical integration methods are implemented in Sage: the midpoint rule, the trapezoidal rule and Simpson's rule. The differences between the exact value of integration and the approximation are tabulated by the number of subintervals $n$ (Fig \ref{fig:integration}).

\small

\begin{verbatim}
f(x) = x^2
a = 0.0
b = 2.0
table = []
exact = integrate(f(x), x, a, b)

for n in [4, 10, 20, 50, 100]:
    h = (b-a)/n
    midpoint = sum([f(a+(i+1/2)*h)*h for i in range(n)])
    trapezoid = h/2*(f(a) + 2*sum([f(a+i*h) for i in range(1,n)])
                + f(b))
    simpson = h/3*(f(a) + sum([4*f(a+i*h) for i in range(1,n,2)])
                + sum([2*f(a+i*h) for i in range (2,n,2)]) + f(b))
    table.append([n, h.n(digits=2), (midpoint-exact).n(digits=6),
       (trapezoid-exact).n(digits=6), (simpson-exact).n(digits=6)])

html.table(table, header=["n", "h", "Midpoint rule",
           "Trapezoidal rule", "Simpson's rule"])
\end{verbatim}

\normalsize

\begin{figure}[h]
	\centering
		\includegraphics[width=0.80\textwidth]{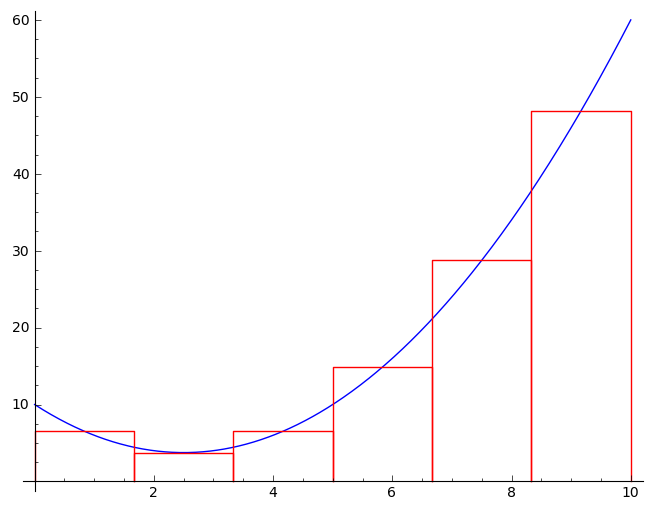}
	\caption{The geometric interpretation of the midpoint rule is visualized using Sage's built-in functions for numerical analysis and plotting.}
	\label{fig:integration2}
\end{figure}

There are also built-in methods for numerical integration in Sage. For instance, it is possible to automatically produce piecewise-defined line functions defined by the trapezoidal rule or the midpoint rule. These functions can be used to visualize different geometric interpretations of the numerical integration methods. In the next example, midpoint rule is used to calculate an approximation for the definite integral of the function $f(x)=x^2-5x+10$ over the interval $[0, 10]$ using six subintervals (Fig \ref{fig:integration2}).

\small
\begin{verbatim}
f(x) = x^2-5*x+10
f = Piecewise([[(0,10), f]])
g = f.riemann_sum(6, mode="midpoint")
F = f.plot(color="blue")
R = add([line([[a,0],[a,f(x=a)],[b,f(x=b)],[b,0]], color="red")
    for (a,b), f in g.list()])
show(F+R)
\end{verbatim}
\normalsize


Sage's numerical integration method \textsl{numerical\_integral} utilizes the adaptive Gauss-Kronrod method available in the GSL \textsl{(GNU Scientific Library)} library. More methods for numerical integration are available in SciPy's sub-packages.

\subsection{Multidimensional Newton's method}

Newton's iteration for solving a system of equations
$f(x) = 0$ in ${\mathbb R}^n$ consists of fixing
 a suitable initial value  $x_0$ and recursively defining $$x_{k+1}=x_k - J_f (x_k)^{-1}f(x_k)\,, k=0,1,2,\dots .$$

Consider next the case $n=3$ and the function
$$f(\mathbf{x})=\left(\begin{array}{c}3x_0-\cos{(x_1 x_2)}-\frac{1}{2}\\x_0^2-81(x_1+0.1)^2+\sin{x_2}+10.6\\e^{-x_0 x_1}+20x_2+\frac{10\pi-3}{3}\end{array}\right),$$
where  $\mathbf{x}=(x_0, x_1, x_2).$ In this program the Jacobian
matrix $J_f(x)$ is computed symbolically and its inverse
numerically. As a result, the program produces a table of the iteration steps and an interactive 3D plot that shows the steps in a coordinate system.


\small
\begin{verbatim}
x0, x1, x2 = var('x0 x1 x2')
f1(x0, x1, x2) = 3*x0 - cos(x1*x2) - (1/2)
f2(x0, x1, x2) = x0^2 - 81*(x1 + 0.1)^2 + sin(x2) + 10.6
f3(x0, x1, x2) = e^(-x0*x1) + 20*x2 + (10*pi - 3)/3
f(x0, x1, x2) = (f1(x0,x1,x2), f2(x0,x1,x2), f3(x0,x1,x2))

j = jacobian(f, [x0,x1,x2])
x = vector([3.0, 4.0, 5.0]) # Initial values
data = [[0, x, n(norm(f(x[0], x[1], x[2])), digits=4)]]

for i in range(1,8):
    x = vector((n(d) for d in x - j(x0=x[0], x1=x[1],
    x2=x[2]).inverse()*f(x[0], x[1], x[2])))
    data.append([i, x, norm(f(x[0], x[1], x[2]))])

# HTML Table
html.table([(data[i][0], data[i][1].n(digits=10),
           n(data[i][2], digits=4)) for i in range(0,8)],
           header = ["$i$", "$(x_0,x_1,x_2)$", "$norm(f)$"])

# 3D Picture
l = line3d([d[1] for d in data], thickness=5)
p = point3d(data[-1][1], size=15, color="red")
show(l + p)
\end{verbatim}
\normalsize

\begin{figure}[ht]
	\centering
		\includegraphics[width=0.90\textwidth]{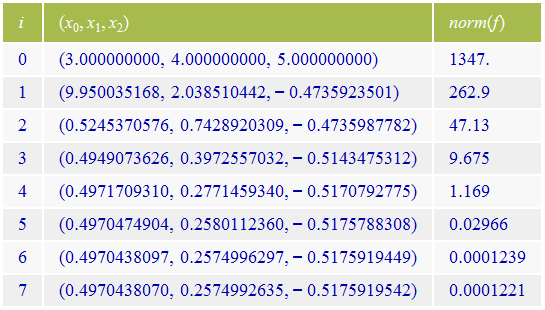}
	\caption{The program produces a table showing the iteration steps of the Newton's method.}
	\label{fig:multinewton1}
\end{figure}


\begin{figure}[ht]
	\centering
		\includegraphics[width=0.70\textwidth]{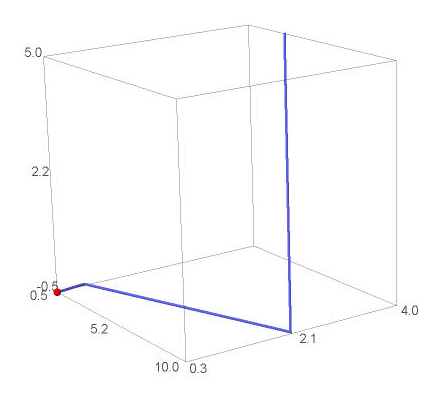}
	\caption{An interactive 3D plot shows the iteration steps in a coordinate system. The plot is made with the Jmol application integrated in Sage.}
	\label{fig:multinewton2}
\end{figure}

\subsection{Nonlinear fitting of multiparameter functions}

Given the data $(x_j, y_j), j=1,\dots,m$, we wish to fit $y=f(x, \lambda)$ into a "model", where $\lambda=(\lambda_1, ..., \lambda_p)$ by minimizing the object function

$$s(\lambda)=\sum_{j=1}^m (y_j-f(x_j, \lambda))^2.$$

The minimization may encounter the usual difficulties:
the minimum need not be unique and there may be several
local minima, each of which could cause the algorithm
to stop prematurely. In the next example the function
\textsl{minimize} uses the Nelder-Mead Method from the
\textsl{scipy.optimize} package.

   Let $\lambda=(\lambda_1, \lambda_2, \lambda_3)$.
   Consider the model function $$f(x, \lambda) = \lambda_1 e^{-x}+\lambda_2 e^{-\lambda_3 x}.$$

The data points used in this example are generated
randomly by deviating the values of the model function.

\clearpage

\small
\begin{verbatim}
from numpy import *

def fmodel(lam, x):
    return lam[0]*exp(-x) + lam[1]*exp(-lam[2]*x)

def fobj(lam, xdata, ydata):
    return linalg.norm(fmodel(lam, xdata) - ydata)

xdata = arange(0, 1.15, 0.05)
lam = [0.2, 1.5, 0.7]
y = fmodel(lam, xdata)

# The generation of the data points
ydata = y*(0.97+0.05*random.rand(y.size))

# Initial values
lam0 = [1, 1, 1]
y0 = fobj(lam0, xdata, ydata)

# The minimization of the object function
lam = minimize(fobj, lam0, args=(xdata, ydata), algorithm='simplex')

yfinal = fobj(lam, xdata, ydata)

# Plot of the datapoints and the model function
fit = plot(fmodel(lam, x), (x, 0, 1.5), legend_label="Fitted curve")
datapoints = list_plot(zip(xdata, ydata), size=20,
                       legend_label="Data points")

html("\n\n$\\text{Object function values: start = %s, final = %s}$\n"
     %( n(y0, digits=5), n(yfinal, digits=5)))
show(fit + datapoints, figsize=5, gridlines=True,
     axes_labels=("xdata", "ydata"))
\end{verbatim}
\normalsize

\begin{figure}[ht]
	\centering
		\includegraphics[width=0.90\textwidth]{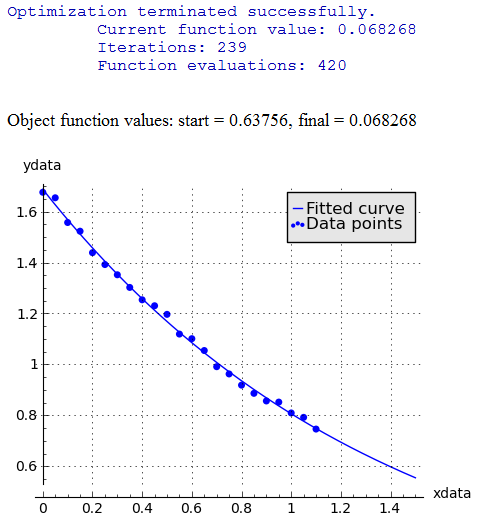}
	\caption{The algorithm used in the program returns a report on the success of the optimization. The plot shows the data points and the model function in the same coordinate system.}
	\label{fig:optimization}
\end{figure}


\subsection{Polynomial Approximation}

The problem of finding $(n-1)$th order polynomial approximation for a function $g$ on the
interval $[r_1,r_2]$ leads to the minimization of the expression $$f(c_1, ..., c_n)=\int_{r_1}^{r_2} (g(x)-\sum_{k=1}^n c_k x^{n-k})^2\ dx$$
with respect to the parameter vector
$(c_1,...,c_n)\,.$
In order to find the optimal value of the parameter vector, we consider the critical points where gradient vanishes
i.e. the points where
 $$\frac{\partial f}{\partial c_i}=0\,, \forall i = 1, ..., n\,.$$

 For the purpose of illustration, consider the case when $g(x)=e^x$ and $n = 2, 3, 4$. The equations $\frac{\partial f}{\partial c_i}=0$ lead to the requirement $$\sum_{k=1}^n c_k \left(\frac{r^{2n-k-j+1}}{2n-k-j+1}\right) \biggr|_{r_1}^{r_2} =\int_{r_1}^{r_2} g(x)x^{n-j}\ dx\,,$$
which is an $n\times n$ linear system of equations for the coefficients
 $c_k\,.$  In the code below the
 integrals on the right hand side are
 evaluated in terms of the function numerical\_integral.

\small
\begin{verbatim}
from numpy import arange, polyval, zeros, linalg
f(x) = e^x
interval = [-1, 1]
nmax = 3
data = zeros((nmax, nmax))
r1, r2 = interval[0], interval[1]

for n in range(2, nmax+1):
    a, b, c = zeros((n, n)), zeros(n), zeros(n)
    for j in range(1, n+1):
        for k in range(1, n+1):
            a[j-1, k-1] = (r2^(2*n-k-j+1) - r1^(2*n-k-j+1))/(2*n-k-j+1)
        b[j-1] = numerical_integral(f*x^(n-j), r1, r2)[0]
    c = linalg.solve(a,b)
    h = (r2-r1)/40
    xs = arange(r1, r2+h, h)
    y1 = [f(xi) for xi in xs]
    y2 = polyval(c, xs)
    err = abs(y2-y1)
    maxer = max(err)

    # Use trapezoidal rule to compute error
    int1 = h*(sum(err) - 0.5*(err[0] + err[-1]))
    int2 = h*(sum(err^2) - 0.5*(err[0]^2 + err[-1]^2))

    # Plots
    eplot = plot(f(x), (x, r1, r2), color="black")
    polyplot = plot(polyval(c, x), (x, r1, r2), color="red", figsize=3)
    epoly = eplot + polyplot
    errplot = plot(abs(f(x)-polyval(c, x)), (x, r1, r2), figsize=3)

    # Output text and graphics
    html("<hr>$n=%s:$"%n)
    html.table([["$%s$"%latex(matrix(a).n(digits=4)),
               "$%s$"%latex(vector(b).column().n(digits=4)),
               "$%s$"%latex(vector(c).column().n(digits=4))]],
               header=["$a$", "$b$", "$c$"])
    html("$\\text{Abs. error = } %s\qquad\qquad\\text{L2 error = }
         %s$"%(maxer, int2))
    html.table([["$\\text{Approximation (n = %s)}$"%n,
               "$\\text{Abs. error function (n = %s)}$"%n],
               [epoly, errplot]], header=True)
\end{verbatim}

\normalsize


\begin{figure}[ht]
	\centering	 \includegraphics[height=0.85\textheight,keepaspectratio]{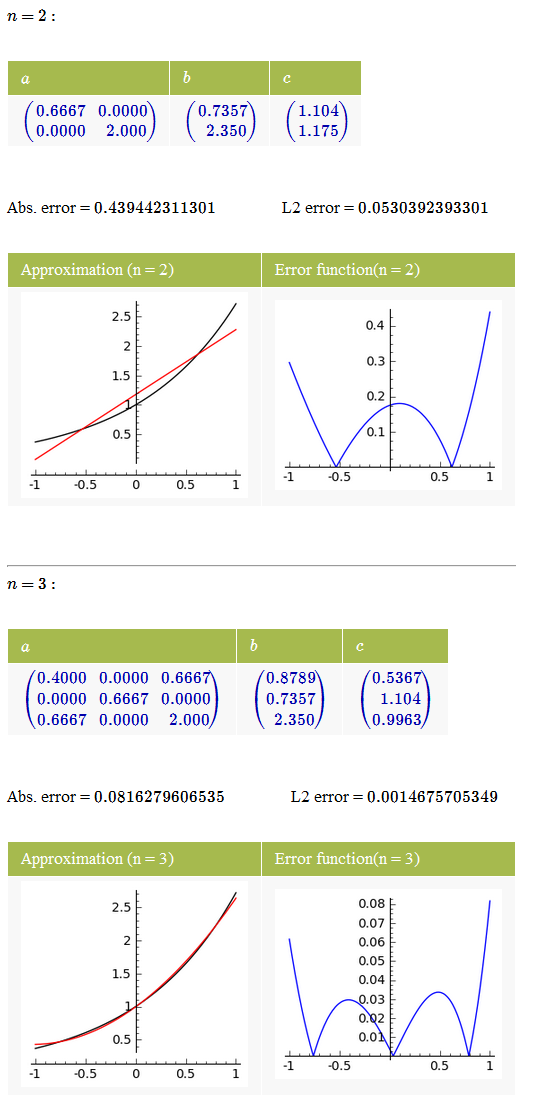}
	\caption{The picture shows the $(n-1)$th order polynomial approximation for the function $e^x$ on the
interval $[-1,1]$ in the cases of $n=2$ and $n=3$.}
	\label{fig:polynomial}
\end{figure}


\section{Concluding remarks}

During its initial years of development, the Sage project has grown to an environment which offers an attractive alternative for the commercial packages in several areas of computational mathematics. For the purpose of scientific computation teaching, the  functionalities of Sage are practically the same as those of commercial packages. While free availability to instructional purposes is a very significant advantage, there
are also other important factors from the learner's point of view:
\begin{itemize}
\item[(1)] The Python language can be used also for many other purposes not tied with the scientific computing. A wide selection of extensions and other special libraries are available in the Internet.
\item[(2)] The support of advanced data structures and support of object-oriented data types and modular program structure is available.
\item[(3)] There is an active users' forum.
\end{itemize}
It is likely that the Sage environment in education will become more popular on all levels of mathematics education from junior high school to graduate level teaching at universities. The support of symbolic computation via Maxima and various numerical packages are noteworthy in this respect. For purposes of teaching scientific computing, the Sage environment and the modules it contains form an excellent option.



\small

\begin{thebibliography}{CdB}

\bibitem[ADV]{ADV} {\sc
M. S. Andersen, J. Dahl, and L. Vandenberghe: }
	CVXOPT: A Python package for convex optimization,
	http://abel.ee.ucla.edu/cvxopt.
	
\bibitem[BBSE]{BBSE} {\sc
R. Bradshaw, S. Behnel, D. S. Seljebotn, G. Ewing, et al.: }
	The Cython compiler,
	http://cython.org.

\bibitem[BF]{BF} {\sc
R. L. Burden and J. D. Faires: } Numerical analysis.  2002

\bibitem[CdB]{CdB} {\sc
S. D. Conte and  C. de Boor:} Elementary numerical analysis: An algorithmic approach. Third ed.  McGraw-Hill Book Co., New York-Toronto, Ont.-London 1965 x+278 pp.

\bibitem[DGPS]{DGPS} {\sc
W. Decker, G.-M. Greuel, G. Pfister and H. Sch{\"o}nemann:}
	Singular --	A computer algebra system for polynomial computations,
	http://www.singular.uni-kl.de.

\bibitem[GAP]{GAP}
	GAP -- Groups, Algorithms, and Programming,
	The GAP Group,
	http://www.gap-system.org.

\bibitem[H]{H} {\sc M. T. Heath:} Scientific computing: An Introductory Survey. Second ed. McGrawHill 2002.

\bibitem[His]{His}
	Sage Reference v5.2: History and License,
	The Sage Development Team,
	2012,
	http://www.sagemath.org/doc/reference/history\_and\_license.html.


\bibitem[Hun]{Hun} {\sc J. D. Hunter:}
	Matplotlib: A 2D Graphics Environment.
	Computing in Science \& Engineering, Vol. 9, No. 3. (2007), pp. 90-95,
	doi:10.1109/MCSE.2007.55.


\bibitem[Jmol]{Jmol}
	Jmol: an open-source Java viewer for chemical structures in 3D,
	http://www.jmol.org.



\bibitem[K]{K} {\sc J.
Kiusalaas: } Numerical methods in engineering with Python. Second edition. Cambridge University Press, New York, 2010.



\bibitem[L]{L}  {\sc H. P.  Langtangen:} A Primer on Scientific Programming With Python, Springer, 2009,
ISBN:	3642024742, 9783642024740

\bibitem[MF]{MF} {\sc
J. H. Mathews and K. D. Fink:} Numerical methods using MATLAB, Third ed. 1999, Prentice Hall, Inc., Englewood Cliffs, NJ.

\bibitem[Mol]{Mol} {\sc
C. B. Moler:} Numerical computing with MATLAB. Society for Industrial and Applied Mathematics, Philadelphia, PA, 2004. xii+336 pp. ISBN: 0-89871-560-1

\bibitem[Num]{Num}
Numerical Computing with Sage,
Release 5.2,
The Sage Development Team,
2012,
http://www.sagemath.org/pdf/numerical\_sage.pdf

\bibitem[O]{O} {\sc
T. E. Oliphant:}
	Python for Scientific Computing,
	Computing in Science \& Engineering 9, 90 (2007).

\bibitem[PARI]{PARI}
	PARI/GP,
	Bordeaux,
	2012,
	http://pari.math.u-bordeaux.fr.

\bibitem[PFTV]{PFTV} {\sc
W. H. Press, S. A. Teukolsky, W. T. Vetterling, and B. P. Flannery:} Numerical recipes. The art of scientific computing. Third edition. Cambridge University Press, Cambridge, 2007. xxii+1235 pp. ISBN: 978-0-521-88068-8

\bibitem[PG]{PG} {\sc F. Pérez and B. E. Granger:}
	IPython: A System for Interactive Scientific Computing,
	Computing in Science \& Engineering 9, 90 (2007).

\bibitem[R]{R}
	R: A Language and Environment for Statistical Computing,
	R Core Team,
	R Foundation for Statistical Computing,
	Vienna, Austria,
	ISBN: 3-900051-07-0,
	http://www.R-project.org.
	
\bibitem[Ras]{Ras} {\sc  A. Rasila:}
  Introduction to numerical methods with Python language, part 1, Mathematics
Newsletter / Ramanujan Mathematical Society 14: 1 and 2 (2004), 1 -15. http://www.ramanujanmathsociety.org/

\bibitem[Sage]{Sage} {\sc
William A. Stein et al.:}
Sage Mathematics Software (Version 5.2),
	The Sage Development Team,
	2012,
	http://www.sagemath.org.

\bibitem[SageTeX]{SageTeX} {\sc
Dan Drake et al.:}
The SageTeX Package,
	2009,
	 ftp://tug.ctan.org/pub/tex-archive/macros/latex/contrib/sagetex/sagetexpackage.pdf.

\bibitem[SciPy]{SciPy} {\sc
E. Jones, T. Oliphant, P. Peterson, et al.:}
	SciPy: Open source scientific tools for Python,
	http://www.scipy.org/.

\bibitem[T]{T} {\sc
S. Tosi:}
Matplotlib for Python Developers,
  From technologies to solutions,
  2009,
  Packt Publishing.

\bibitem[Tac]{Tac} {\sc J. E. Stone:}
	The Tachyon 3D Ray Tracer,
	Sage Reference v5.2,
	The Sage Development Team,
	 http://www.sagemath.org/doc/reference/sage/plot/plot3d/tachyon.html.
	
\bibitem[TLN]{TLN}  {\sc A.
Tveito, H. P.  Langtangen,  B. F. Nielsen, and X. Cai:}  Elements of scientific computing. Texts in Computational Science and Engineering, 7. Springer-Verlag, Berlin, 2010. xii+459 pp. ISBN: 978-3-642-11298-0



\end{thebibliography}
\end{document}